%% file: jsc-3474.tex
\documentclass{elsart}
\usepackage{amssymb,natbib,amsmath,enumerate,latexsym,array}

\usepackage{color}

\usepackage[dvips]{graphicx}
\usepackage{xypic}
\xyoption{all}

\newcommand{\id}   {\mathsf{id}}
\newcommand{\init} {\mathsf{init}}
\newcommand{\norm} {\mathsf{norm}}
\newcommand{\sink} {\mathsf{sink}}

\newcommand{\ep} {\varepsilon}
\newcommand{\bA} {\mathsf{A}}
\newcommand{\ubA}{\underline{\bA}}
\newcommand{\ubAG}{\underline{\bA}_{\bG}}
\newcommand{\ubAl}{\underline{\bA}_{\lambda}}
\newcommand{\ubAT}{\underline{\bA}_{\calT}}
\newcommand{\bG} {\mathsf{G}}
\newcommand{\bH} {\mathsf{H}}
\newcommand{\bK} {\mathsf{K}}
\newcommand{\bL} {\mathsf{L}}
\newcommand{\bM} {\mathsf{M}}
\newcommand{\bN} {\mathsf{N}}
\newcommand{\ubN}{\underline{\bN}}
\newcommand{\ubNG}{\underline{\bN}_{\bG}}
\newcommand{\bT} {\mathsf{T}}
\newcommand{\bTp}{\mathsf{T}_+}
\newcommand{\sets}{\mathsf{Sets}}

\newcommand{\lan}{\langle}
\newcommand{\ran}{\rangle}

\newcommand{\src}{\mathsf{src}}
\newcommand{\tgt}{\mathsf{tgt}}

\newcommand{\RbG}{R_{\bG}}
\newcommand{\RbH}{R_{\bH}}
\newcommand{\RbK}{R_{\bK}}
\newcommand{\RbHK}{R_{\bH\bK}}
\newcommand{\SbG}{S_{\bG}}
\newcommand{\SbH}{S_{\bH}}
\newcommand{\SbK}{S_{\bK}}
\newcommand{\SbHK}{S_{\bH\bK}}
\newcommand{\XbG}{X_{\bG}}
\newcommand{\XbH}{X_{\bH}}
\newcommand{\XbK}{X_{\bK}}
\newcommand{\YbH}{Y_{\bH}}
\newcommand{\YbK}{Y_{\bK}}
\newcommand{\HGK}{\bH\backslash\bG/\bK}
\newcommand{\HGammaK}{\bH\backslash\Gamma/\bK}

\newcommand{\mon}{\mathrm{mon}}
\newcommand{\ob} {\mathrm{Ob}}

\newcommand{\calT}{\mathcal T}

\newcommand{\sfl}{\mathsf{l}}
\newcommand{\pl} {\mathsf{pl}}
\newcommand{\ppl}{\mathsf{ppl}}
\newcommand{\sul} {\mathsf{sl}}
\newcommand{\psl}{\mathsf{psl}}
\newcommand{\irrRCT}{\mathsf{irr}_{R^C}(\bT)}

\newcommand{\circto}{\;\circ\!\!\!\!\!\!\to\!}

\def\leq{\leqslant}

\newtheorem{example}{Example}[section]
\newtheorem{Def}[example]{Definition}

\newenvironment{proof}{\noindent {\bf Proof} }{ \hfill
$\Box$ \mbox{}}

\begin{document}

\begin{frontmatter}

\title{String rewriting for Double Coset Systems}

\author[Bangor]{Ronald Brown\thanksref{label1}}
\ead{r.brown@bangor.ac.uk},
\author[Leicester]{Neil Ghani}
\ead{n.ghani@mcs.le.ac.uk},
\author[Leicester]{Anne Heyworth\thanksref{label2}}
\ead{a.heyworth@mcs.le.ac.uk},
\author[Bangor]{Christopher D. Wensley}
\ead{c.d.wensley@bangor.ac.uk}
\thanks[label1]{INTAS
Project 94-436 ext `Algebraic K-theory, groups and categories'.
{Brown was supported for part of this research by a
Leverhulme Emeritus Fellowship 2002-2004.} }
\thanks[label2]{EPSRC GR/R29604/01 `Kan - A Categorical
Approach to Computer Algebra'.}
\address[Bangor]{Mathematics Division, Informatics, University of Wales,
Bangor, LL57 1UT.}
\address[Leicester]{Department of Computer Science, Leicester University,
Leicester,
LE1 7RH.}

\title{}

\begin{abstract}
In this paper we show how string rewriting methods can be applied to
give a new method of computing double cosets. Previous methods for
double cosets were enumerative and thus restricted to finite
examples. Our rewriting methods do not suffer this restriction and
we present some examples of infinite double coset systems which can
now easily be solved using our approach. Even when both enumerative
and rewriting techniques are present, our rewriting methods will be
competitive because they i) do not require the preliminary
calculation of cosets; and ii) as with single coset problems, there
are many examples for which rewriting is more effective than
enumeration.

Automata provide the means for identifying expressions for normal
forms in infinite situations and we show how they may be constructed
in this setting. Further, related results on logged string rewriting
for monoid presentations are exploited to show how witnesses for the
computations can be provided and how information about the subgroups
and the relations between them can be extracted. Finally, we discuss
how the double coset problem is a special case of the problem of
computing induced actions of categories which demonstrates that our
rewriting methods are applicable to a much wider class of problems
than just the double coset problem.
\end{abstract}


\end{frontmatter}


\section{Introduction} \label{sec-intro}


Given a group $\bG$ and two subgroups $\bH$ and $\bK$, the double
cosets are the equivalence classes of the relation $\sim$ where $g
\sim g' \Leftrightarrow hgk=g'$ for some $h \in \bH, k \in \bK$.
The set of double cosets is written $\HGK$. Combinatorially, the
double cosets are the orbits of the left action of $\bH$ on the
right cosets $\bG/\bK$, and also the orbits of the right action of
$\bK$ on the left cosets $\bH\backslash\bG$. Double coset techniques
give examples of deep and wide applications of group theoretic
methods in chemistry and physics~\citep{RuchKleine}. For example,
there are applications through Polya's theory of counting, to
considerations of deuterons colliding in scattering theory
\citep{pletsch}. Real semisimple symmetric spaces are often
characterised by a pair of commuting involutions of a reductive
group and many of their properties are studied in this setting --
in this case double cosets are of importance for representation
theory of $p$-adic symmetric $K$-varieties~\citep{loek}. Double coset
computation can be seen as a way of constructing finite quotients
of HNN-extensions of known groups or as a way of constructing
groups given by symmetric presentations~\citep{curtis}.

There are a number of different questions which arise if we
want to compute with double cosets, for example: the enumeration of
the double cosets; finding a set of representatives for them;
deciding questions such as whether a pair of group elements lie
within the same double coset or not; and proving either case. 
There are number of algorithms for computing such double coset problems
but these are regarded as incomplete. 
Indeed, in Section 4.6.8 of the recently published survey book on 
computational group theory~\citep{HEO} we find
\emph{``Unfortunately, no really satisfactory algorithm for solving this
problem has been found to date.''} 
In 1981 the first algorithmic methods for computing double cosets 
were published~\citep{Butler}, and applied to groups of order $\leq 10^4$.  
The approach was a variation on Dimino's algorithm~\citep{dimino} 
for computing a list of elements of a small group.  
In~\citep{HEO}, the basic approach for permutation groups 
is to use orbit methods to compute left or right cosets, 
and then orbits of these cosets to obtain double cosets. 
As pointed out in~\citep{HEO}, this may involve the calculation 
of a large number of cosets in order to determine 
a small number of double cosets.  
More recent methods for computing double cosets of finitely presented groups 
use Todd-Coxeter procedures~\citep{Linton}. 
All of these methods have been implemented in the commonly used programs 
for computational discrete algebra, 
{\sf GAP}~\citep{GAP} and {\sf MAGMA}~\citep{MAGMA}. 

The primary alternative to Todd-Coxeter procedures for
ordinary coset enumeration and computation of groups given by
presentations is string rewriting \citep{paper2, Sims}. 
In finite settings the two approaches are comparable: 
certain problems being more effectively addressed by the enumerative method 
and others benefiting more from a rewriting approach. 
However, for cases involving infinitely many elements, 
rewriting rather than enumeration is the natural choice. 

This paper {demonstrates} how string rewriting can be 
applied to the problems of computing double cosets, giving a new
alternative to the Todd-Coxeter methods {and one which }
can further be applied to infinite groups. In particular, this
paper makes the following contributions.
\begin{itemize}
\item 
The introduction of the notion of a double coset rewriting
system and associated Knuth-Bendix completion algorithm as a
mechanism for attempting to decide whether two elements of a group
lie in the same double coset.
\item 
The specification, in Section \ref{sec-automata}, of a process
which takes a finite complete double coset rewriting system and
constructs a finite state automaton whose language is a set of
unique normal forms for the double cosets. 
\item 
The specification of a higher dimensional version of the 
Knuth-Bendix algorithm and logged double coset rewriting. 
This gives, for example,
presentations for the subgroups defining the double cosets. 
\item 
A discussion of the implementation of these algorithms as a
deposited package {\sf kan} for {\sf GAP}. 
\item
In Section \ref{sec-kan} we put our algorithm in context by showing
how it arises as a special case of rewriting for an induced action 
of categories, using a category $\bH\circto\bK$ constructed from 
the two subgroups.
\end{itemize}

The authors would like to thank Steve Linton for helpful discussions,
especially regarding the {\sf GAP} implementations of the methods
presented in this paper; 
and Tim Porter for help with applications of category theory.
Thanks are also due to the referee for several helpful suggestions.  
Figures were typeset using {\sf XFig}.


\section{Rewriting for Double Cosets} \label{sec-dc}

If we consider using rewriting to solve double coset problems
we have a choice: to develop a specialised type of rewriting
for this situation, or to rephrase the problem in a way that
allows existing techniques to be applied.
The former method has the advantage of specialty --
the ``double coset rewriting systems'' can be examined in isolation
as though they were in some way an advance on existing methods
rather than a useful application.
This is certainly of some value if one is wishing to write a very
specialised program, designed to compute only double coset problems and
investigate the particular properties of rewriting systems of this type,
but it can obscure the simplicity and the best features of the result.
Therefore, we choose to adopt the most straightforward method --
simulating the required computations by embedding the group $\bG$
in a particular free monoid and then applying standard procedures~\citep{Book}.
We then have to show that the structure we wish to compute
coincides with the rewriting model used.

\begin{Def}[Presentation of a Double Coset System]{\mbox{ }}\\
  Let $\bG$ be a group with monoid presentation $\mon\langle \XbG,
  \RbG \rangle$ and let $\theta: \XbG^* \to \bG$ be the natural monoid
  homomorphism.  Let $\XbH, \XbK \subseteq \XbG^*$ be such that
  $\YbH=\theta(\XbH)$ and $\YbK=\theta(\XbK)$ are sets of generators
  for the subgroups $\bH$ and $\bK$ respectively.  Then we shall say
  that the quadruple $(\XbG, \RbG, \XbH, \XbK)$ is a presentation of
  the system of double cosets $\HGK$.
\end{Def}

If $R$ generates an equivalence relation or congruence on a free
monoid $S$ then the class of $s \in S$ is denoted $[s]_R$.  
Similarly, we write $\HGK = \{[g]_{\sim} ~|~ g \in G\}$, 
where $\sim$ is the relation defined at the beginning of Section 1.  
If $\mon\lan \XbG, \RbG \ran$ is a monoid presentation for $\bG$ 
then the free monoid $\bTp$ in which we compute 
is generated by $\XbG$ together with two extra (tag) variables $H$ and $K$.
A string $HwK$ represents the double coset $[\theta w]_{\sim}$,
and we require a congruence $\stackrel{*}{\leftrightarrow}_R$
such that $HwK \stackrel{*}{\leftrightarrow}_R Hw'K$
if and only if $\theta w' \sim \theta w$.
Clearly $\bTp$ contains many elements that are not of the form $HwK$,
but these do not arise in the computations performed
when completing the rewriting system which determines the double cosets.
This is the key observation which allows us to use standard methods.

\begin{thm}[Double Coset Rewriting]\label{main}{\mbox{ }}\\
Let $(\XbG, \RbG, \XbH, \XbK)$
be the data for the double coset system $\HGK$,
let $H$ and $K$ be symbols, and let $\bT$
be the subset of terms of the form $H w K$ of the free monoid
$\bTp =(\{H,K\} \cup \XbG)^*$ where $w \in \XbG^*$.
Define
$$
R = \RbG \cup \{(H h, H): h \in \XbH\} \cup \{(k K, K): k \in \XbK\}.
$$
Let $\to_R$ be the reduction relation generated by $R$ on the free monoid
$\bTp$, and let $\stackrel{*}{\leftrightarrow}_R$
be the reflexive, symmetric, transitive closure of $\to_R$,
which is the congruence generated by $R$.
Then there is a bijection of sets
$$
\frac{\bT}{\stackrel{*}{\leftrightarrow}_R}
\quad\cong\quad
\frac{\bG}{\sim}.
$$
\end{thm}

\begin{proof}
We first show that there is a well-defined map
$$
\phi \;:\;
\frac{\bT}{\stackrel{*}{\leftrightarrow}_R}
\to \frac{\bG}{\sim} \;
\quad\text{ where }\quad
\phi([H w K]_R) = [\theta w]_\sim \;.
$$
If $[H w K]_R = [H w' K]_R$, there exists a sequence $w=w_1,w_2
\ldots, w_n=w'$ in $\XbG^*$ such that for $i=1, \ldots, n\!-\!1$
either $(w_i,w_{i+1})$ or $(w_{i+1},w_i)$ 
has one of the following forms:
\begin{enumerate}[i)]
\item
$(ulv,urv)$ for some $(l,r) \in \RbG$, $u, v \in \XbG^*$,
\item
$(hv,v)$ for some $h \in \XbH$, $v \in \XbG^*$,
\item
$(uk, u)$ for some $k \in \XbK$, $u \in \XbG^*$.
\end{enumerate}
Since $\theta$ is a monoid homomorphism, in the first case
$\theta(w_i)=\theta(w_{i+1})$, in the second case $\theta(hv) \sim
\theta(v)$ and in the third case $\theta(uk) \sim \theta(u)$.  Thus in
all cases $[\theta(w_i)]_\sim = [\theta(w_{i+1})]_\sim$ as required.

Secondly, let $\tau : \bG \to \XbG^*$ be a section of $\theta$ 
so that $(\theta \circ \tau)g = g$ 
and $(\tau \circ \theta)g \stackrel{*}{\leftrightarrow}_{\RbG} g$ for all $g \in \bG$. 
Further,
$\tau(g_1g_2) \stackrel{*}{\leftrightarrow}_{\RbG} \tau(g_1)\tau(g_2)$
since $\theta$ maps both $\tau(g_1)\tau(g_2)$ and $\tau(g_1g_2)$ 
to $g_1g_2$.
Define 
$$
\phi' \;:\; \frac{\bG}{\sim} \to
\frac{\bT}{\stackrel{*}{\leftrightarrow}_R} \quad\text{ where
}\quad \phi'([g]_\sim) = [H (\tau g) K]_R.  $$
To verify that $\phi'$
is well-defined, suppose $[g]_\sim = [g']_\sim$ for some $g, g' \in
\bG$.  Then, by the definition of $\sim$, we have $h \in \XbH^*$ and
$k \in \XbK^*$ such that $g' = (\theta h) g (\theta k)$, so that 
$$
H (\tau g') K \stackrel{*}{\leftrightarrow}_R 
H (\tau \theta h)(\tau g)(\tau \theta k) K
\stackrel{*}{\leftrightarrow}_{\RbG}~ H h (\tau g) k K
\stackrel{*}{\to}_R~ H (\tau g) K.  
$$
Finally, we observe that
$\phi( \phi'([g]_\sim ))
= \phi( [H (\tau g) K)]_R )
= [\theta(\tau g)]_\sim = [g]_\sim$,
and that
$\phi'( \phi( [H w K]_R))
= \phi( [(\theta w)]_\sim)
= [ H (\tau (\theta w)) K]_R
= [ H w K]_R$
since $\RbG \subseteq R$.
Thus $\phi$ is a bijection with inverse $\phi'$.
\end{proof}

Note the use of the tags $H$ and $K$. They provide a particularly
simple way of deleting elements of $\XbH$ from a word providing they
occur at the far left of the word and similarly for deleting elements
of $\XbK$ providing they occur at the far right of the word.  
Given the double coset presentation $(\XbG, \RbG, \XbH, \XbK)$ 
as in the theorem above, we may refer to $R$ as a 
\emph{double coset rewriting system} for $\HGK$.  
Of course $R$ may not be complete (confluent and noetherian) 
and so the natural next step would be to use Knuth-Bendix completion 
to try and obtain an equivalent, but complete, rewriting system. 
As justified above, we choose to perform this completion 
by considering rewriting over the free monoid $\bTp$ 
rather than the set $\bT$.  However, we must then be
sure that if $\to_{R'}$ is complete on $\bTp$, and its
closure $\stackrel{*}{\leftrightarrow}_{R'}$ coincides with
$\stackrel{*}{\leftrightarrow}_{R}$, then the restriction of
$\to_{R'}$ to $\bT$ is also complete.  Fortunately, this is
obviously true as $\to_R$ is closed on the subset $\bT$, that
is, if $w_1 \to_R w_2$ then $w_1 \in \bT$ if and only if $w_2 \in
\bT$. To see this, note that $\leftrightarrow_{R}$ never
removes or adds tags. Thus, we can immediately state the following
corollaries.


\begin{cor}
If the double coset rewriting system $R$ can be completed on $\bTp$
then we have a solution to the problem of deciding whether two elements
$g, g'$ of $\bG$ lie within the same double coset.
\end{cor}

\begin{cor}
If the double coset rewriting system $R$ can be completed on
$\bTp$
then we can find a unique normal form for each double coset.
\end{cor}

\begin{rem}[Implementations]\mbox{ }\\
\emph{Besides the fact that the string rewriting methods
we have presented enable us to tackle problems involving
infinite groups, they also allow us to immediately use
existing string rewriting programs such as
those in {\sf GAP} and in {\sf KBMAG}~\citep{KBMAG} to compute double cosets.}

\emph{An alternative approach is to complete a rewriting system for $\bG$ 
and then construct a double coset rewriting system on the elements of $\bG$.  
Such a system obeys more subtle laws than a standard rewriting system 
on a free monoid.  
For example, if $g>g'$ in the termination order 
used by the completion process, we may not deduce that $gg^{-1}$ 
is greater than $g'g^{-1}$.  
Consequently, we believe the approach we have chosen is cleaner 
than if we were to have worked with multiple rewriting systems 
at different levels to describe the one structure.}

\emph{While it is straightforward to try to use a standard 
Knuth-Bendix completion program to calculate the completion of a 
double coset rewriting relation on the free monoid $\bT_+$, 
it will necessarily be less efficient than a specialised Knuth-Bendix 
program which, for example, restricts itself to the non-free monoid $\bT$.  
For example there will be tests for more overlaps between words 
$wK$ and $zK$ than can possibly arise: we know (but the program does not) 
that the only way in which an overlap can occur is when $z=uw$ or $w=uz$, 
since the tag symbol $K$ will not occur within the strings $w$ or $z$.
A specialised program would also allow different types of ordering,
treating symbols $H,K$ in a way different from those in $\XbG$,
yielding results that could not be obtained with a standard ordering. 
If one wishes to do many calculations of this type, 
it would be worth refining the system to recognise tags and deal with
tagged rules sensibly, and to separate rules into subsystems which
are completed separately. 
Such an approach would not be designed specifically for double cosets 
and could have many other applications~\citep{paper2}.}
\end{rem}


\section{Completion of Rewriting Systems for Double Cosets}
\label{sec-complete}

As we have seen, if the double coset rewriting system $R$ is complete,
we can solve problems such as whether two elements of the group belong to
the same double coset. 
Usually, $R$ is incomplete, and so we attempt to convert it 
to an equivalent complete system.  
We can apply the Knuth-Bendix completion procedure~\citep{K-B} 
to $R$ on $\bTp$ in the standard way, as detailed below.
If $R$ completes on $\bTp$, we are required to prove
that the restriction of $\to_R$ to $\bT$
is preserved throughout the algorithm.

\begin{alg}[Completion]{\mbox{ }}
\vspace{-1ex}
\begin{enumerate}[{\rm K1}]
\item
{\bf (Input)~}
Start with a set of pairs
$R = \RbG \cup \{(H h, H) : h \in \XbH\} \cup \{(k K,K)
   : k \in \XbK\} \subset \bTp$
and a compatible well-ordering on $\bTp$.
\item {\bf (Search)} Find all overlaps: pairs of rules $(l_1,r_1)$,
  $(l_2,r_2)$ which may be applied to the same word and which coincide
  on some subword.  There are essentially two cases, $u l_1 v = l_2$
  or $u l_1 = l_2 v$ for some $u, v \in \bT$.  Add each pair
  $(u r_1 v, r_2)$ or $(u r_1, r_2 v)$ to a list of critical pairs.
\item
{\bf (Resolve)~}
Attempt to resolve each critical pair by reducing both terms with respect
to the rules in $R$.
If the reduced terms are equal then the pair has resolved,
otherwise the reduced pair is orientated according to the well-ordering
and added to a set of new rules and to $R$.
\item
{\bf (Loop)~}
If no new rules were added then go to the next step.
Otherwise, repeat the last two steps with the new set of rules,
checking pairs that arise between the new rules
and between the new rules and the old rules,
but not pairs just between old rules (as these have already been checked).
\item
{\bf (Output)~}
Return the resulting $R^C$, a complete rewriting system on $\bTp$
with respect to the given well-ordering.
\end{enumerate}
\end{alg}

We now prove that if the input for the algorithm is a
double coset rewriting system of the form specified earlier,
then the algorithm will attempt to produce an equivalent complete system.

\begin{thm}[Completeness of Double Coset Rewriting Systems]{\mbox{ }}\\
Let the input for the above algorithm be a double coset
rewriting system $R$ as given in Theorem \ref{main},
such that the algorithm terminates, giving output $R^C$.
Then the restriction of $\to_{R^C}$ to $\bT$
is a complete rewriting system equivalent to the restriction of
$\to_R$ to $\bT$.
\end{thm}

\begin{proof}
We prove the result directly, by showing that no step in the
completion procedure alters the restriction of
$\stackrel{*}{\leftrightarrow}_R$ to $\bT$.
The argument holds because
of the form of the input and the way in which new pairs are generated.

It is convenient (at each stage of the algorithm) to partition a
rewrite system $R$ into four subsets $\RbH$, $\RbK$, $\RbG$ and
$\RbHK$, depending on whether the rule involves $H$, $K$, neither or
both.  The subset $\RbHK$ is initially empty.  Formally:
\begin{eqnarray*}
\RbG &~:=~& \{ (l, r) \in R \mid l, r \in \XbG^* \}, \\
\RbH &~:=~& \{ (Hl, Hr) \in R \mid l, r \in \XbG^* \}, \\
\RbK &~:=~& \{ (lK, rK) \in R \mid l, r \in \XbG^* \}, \\
\RbHK &~:=~& \{ (HlK, HrK) \in R \mid l, r \in \XbG^* \}.
\end{eqnarray*}
We observe that in step K3 of the algorithm a new rule is generated from an
overlap of the left hand sides of two existing rules
(followed by their subsequent reductions). \\

Overlaps of rules $(l_1,r_1)$ and $(l_2,r_2)$ in $\RbG$
may be separated into five types.
In Table 1,  
$l_1$ is either a prefix of $l_2$;
a suffix of $l_2$; an internal subword of $l_2$;
or the overlap is offset to the left or the right.
Neither of $u,v$ may equal the empty word $\id$.

\begin{table}[htb] \label{table-overlapsG}
\begin{center}
\begin{tabular}{|r|c|c|c|}
\hline
type & overlap case & new rule to add & picture \\
\hline
prefix &
$l_1 v = l_2$ &
$(r_1 v, r_2)$ &
$\stackrel{\text{---}~~~~}{\text{---}\!\!\text{---}~}$ \\
suffix &
$u l_1 = l_2$ &
$(u r_1, r_2)$ &
$\stackrel{~~~~\text{---}}{~\text{---}\!\!\text{---}}$ \\
internal subword &
$u l_1 v = l_2$ &
$(u r_1 v, r_2)$ &
$\stackrel{~~\text{---}~~}{\text{---}\!\!\text{---}\!\!\text{---}}$ \\
left offset &
$l_1 v = u l_2$ &
$(r_1 v, u r_2)$ &
$\stackrel{\text{---}\!\!\text{---}~~}{~~\text{---}\!\!\text{---}}$ \\
right offset &
$u l_1 =  l_2 v$ &
$(u r_1, r_2 v)$ &
$\stackrel{~~\text{---}\!\!\text{---}}{\text{---}\!\!\text{---}~~}$ \\
\hline
\end{tabular}

\medskip
Table 1 ~:~ Five types of overlap in $\bG$.
\end{center}
\end{table}

Using the same terminology for the other types of rule,
we see that the overlaps which may occur are of the types shown in
Table 2.  

\begin{table}[hbt]  \label{table-overlapsGHK}
\begin{center}
\begin{tabular}{|c|c|c|c|c|c|}
\hline
overlap words & prefix & suffix & internal & left offset & right offset \\
\hline
$ l_1 ,  l_2 $ 
  & $\checkmark$ & $\checkmark$ & $\checkmark$ & $\checkmark$ & $\checkmark$ \\
$ l_1 , Hl_2 $ 
  &              & $\checkmark$ & $\checkmark$ &     & $\checkmark$ \\
$ l_1 ,  l_2K$ 
  & $\checkmark$ &              & $\checkmark$ & $\checkmark$ &     \\
$Hl_1 , Hl_2 $ 
  & $\checkmark$ &              &              &              &     \\
$Hl_1 ,  l_2K$ 
  &              &              &              & $\checkmark$ &     \\
$ l_1K,  l_2K$ 
  &              & $\checkmark$ &              &              &     \\
$ l_1 , Hl_2K$ 
  &              &              & $\checkmark$ &              &     \\
$Hl_1 , Hl_2K$ 
  & $\checkmark$ &              &              &              &     \\
$ l_1K, Hl_2K$ 
  &              & $\checkmark$ &              &              &     \\
\hline
\end{tabular}

\medskip
Table 2 ~:~ Overlaps types for all pairs of rules.

\medskip
\end{center}
\end{table}

In each case we observe that the new pair added to $R$
will not change the definition of $\stackrel{*}{\leftrightarrow}_R$
when it is restricted to $\bT$
since the new pair is an element of $\stackrel{*}{\leftrightarrow}_R$
and also an element of $\bT \times \bT$.

The reduction of the new rules with respect to the existing rules
also gives a pair which is already an element of
$\stackrel{*}{\leftrightarrow}_R$ when it is restricted to $\bT$.
This is because replacement of a substring $l_i$ in a word in $\bTp$
by $r_i$ will preserve the positions of the tags.
\end{proof}

\section{Automata for Double Coset Rewriting Systems} \label{sec-automata}

When the number of double cosets is infinite,
it is not possible to list all their normal forms,
so a regular expression giving an impression
of the shape or pattern of these forms may be useful.
In simpler string rewriting systems, for example with left or right
cosets, we can build up a catalogue of the normal forms.
In finite cases where we have a length non-increasing order
this is extremely effective, and even in infinite cases
it can serve to show up patterns.
The cataloguing procedure relies heavily on the fact that,
if a term is reducible and generators are appended onto a chosen end,
then the term remains reducible.

In the case of double cosets and doubly tagged strings
we cannot use these methods.
If a term $HlK$ is reducible,
then one of $l \to r$ or $Hl \to Hr$ or $lK \to rK$
or $HlK \to HrK$ is true,
but we cannot assume that either $HlxK$ or $HxlK$ are reducible.
For example, if $Hl \to Hr$, then we are unable to deduce that $Hxl \to Hxr$.
Since we cannot use cataloguing, we turn to the alternative:
in string rewriting we use automata when we wish to find an expression
for the set of normal forms and cataloguing is insufficient.
Here, we describe  techniques for constructing automata
whose languages are the sets of  normal forms for our double cosets.
First we recall the construction of an automaton
for accepting normal forms in $\bG$.

We will use the following notation.
For any set of rules $R$ we set $\sfl(R) = \{l ~|~ (l,r) \in R\}$,
the set of \emph{left-hand sides} of these rules.
Then $\pl(R)$ is the set of all \emph{prefixes} of the rules
and $\ppl(R)$ is the set of all \emph{proper prefixes}:
$$
 \pl(R) = \{u \mid (uv,r) \in R,~ u \not= \id \}, \quad
\ppl(R) = \{u \mid (uv,r) \in R,~ u,v \not= \id \}.
$$
Similarly, $\sul(R)$ and $\psl(R)$ denote the sets of \emph{suffices}
and \emph{proper suffices}.
A non-deterministic automaton $\bN$, 
with state set $S$; alphabet $\Sigma$; initial state $s_0 \in S$; 
transition function $\delta : S \times \Sigma \to 2^S$; 
and accepting states $A \subseteq S$, 
is written $\bN = (S,\Sigma,s_0,\delta,A)$.
A deterministic automaton has $\delta : S \times \Sigma \to S$. 

\begin{Def}[Word acceptor for $\bG$]   \label{def-wag}
\mbox{ }\\
\emph{The \emph{word acceptor} of a group
$\bG = \mon\langle \XbG, \RbG \rangle$
with a finite complete rewrite system $\RbG^C$ is constructed as follows.
First form a non-deterministic automaton
$$
\ubNG ~=~
(~\ppl(\RbG^C) \cup \{\id,\sink\},~ \XbG,~ \id,~ \delta_{\bG},~ \{\sink\})
$$
whose states consist of all proper prefixes $p$ from $\RbG^C$;
the identity word as initial state;
and a sink state $\sink$ which is the only accepting state.}

\emph{The transition function is given by
\begin{eqnarray*}
\delta_{\bG}(\sink,x)
  &~=~&  \{\sink\} \qquad\text{ for all }~ x \in \XbG, \\
\delta_{\bG}(p,x)
  &~=~&  \left\{ \begin{array}{l}
        \{\sink\} \quad\text{ if }~ px=ul
            ~\text{ for some }~ l \in \sfl(\RbG),~
        \text{ else } \\
        \{p_i \in \ppl(\RbG) ~|~ px=u_ip_i
            ~\text{ for some }~ u_i \in \XbG^* ~\}.
        \end{array} \right.
\end{eqnarray*}
Standard results of automata theory~\citep{cohen,lawson}
allow us to determinize $\ubNG$ (using the accessible subset construction); 
take the complement (accepting states become non-accepting, and conversely);
and minimize, giving a deterministic automaton $\ubAG$
which accepts only the normal forms of elements of $\bG$.}
\end{Def}

\begin{thm}[Word acceptor for $\HGK$]\label{thm-dcwa}{\mbox{ }}\\
Let $R^C$ be a finite, complete, double coset rewriting system for subgroups 
$\bH$ generated by $\YbH = \theta\XbH$ 
and $\bK$ generated by $\YbK = \theta\XbK$ 
of the group $\bG$ which is given by the monoid presentation
$\mon\langle \XbG, \RbG \rangle$.
Let $\bT$ and $\bT_+$ be defined as in Theorem \ref{main}.
Then there is a regular expression
representing a regular language $\bL$ over $\bT_+$ such that
$\bL= \{ [g]_\sim : g \in \bG \}$.
\end{thm}

\begin{proof}
We define a non-deterministic automaton $\ubN$
with input alphabet $\Sigma = \XbG \cup \{H,K\}$
which accepts exactly the set $\irrRCT$ of irreducible elements
of $\bT$ with respect to $\to_{R^C}$.

As before, we partition the rules in $R^C$ into
$\RbG^C \cup \RbH^C \cup \RbK^C \cup \RbHK^C$,
where $\RbG^C$ is the complete rewrite system for $\bG$.
The automaton has four main components:
\begin{itemize}
\item
the non-deterministic form $\ubNG$ of the word acceptor for $\bG$,
with states $\SbG$, as in Definition \ref{def-wag}~;
\item
an \emph{$H$-tree}, whose states are $\SbH = \{H\id\} \cup \ppl(\RbH^C)$~;
\item
a \emph{$K$-tree}, whose states are $\SbK = \{\id K\} \cup \psl(\RbK^C)$~;
\item
an \emph{$HK$-tree}, whose states are
$\SbHK = \{H\id \!\cdot\! K\} \cup (\ppl(\RbHK^C) \!\cdot\! K)$~.
\end{itemize}

\medskip
\begin{figure}[hbt] \label{fig-sketch}
\begin{center}
\input{jsc-3474-fig1.pstex_t}

\medskip
Figure 1 ~:~ Sketch of a double coset automaton.
\end{center}

\medskip
\end{figure}

There are two additional states, an initial state $\init$ and a
normal form state $\norm$ which is the only non-accepting state.
An informal sketch showing how these components and states
are connected by transitions is shown Figure 1.

The formal definition of the non-deterministic automaton $\ubN$ is

\qquad
$\ubN ~=~
  (~ S,~ \Sigma = \XbG \cup \{H,K\},~ \init,~ \delta,~ S\setminus\{\norm\})$,
\quad where

\qquad
$S ~=~ \{\init,\norm\} \cup \SbG \cup \SbH \cup \SbK \cup \SbHK$.

The transition function $\delta$ is defined in Table 3,
where $s \in S, x \in \XbG,~ a \in \Sigma,~
p \in \ppl(\RbG^C),~ Hp \in \ppl(\RbH^C),~ qK \in \psl(\RbK^C)$
and $Hp \!\cdot\! K \in \ppl(\RbHK)\!\cdot\! K$.

\medskip
\begin{table}[h]  \label{table-delta}
\begin{center}
\begin{tabular}{|r|rcl|}
\hline
{\bf location} & & & {\bf transition} \\
\hline
\hline
from $\init$
 & $\delta(\init,H)$
   & $=$ & $\{\id, H\id, H\id \!\cdot\! K\}$ \\
 & $\delta(\init,a)$
   & $=$ & $\{\sink\} \quad\text{when}~ a \neq H$ \\
\hline
by $H$
 & $\delta(s,H)$
   & $=$ & $\{\sink\} \quad \text{when}~ s \neq \init$ \\
\hline
by $K$
 & $\delta(p,K)$
   & $=$ & $\{\norm\}$ \\
 & $\delta(\id K,K)$
   & $=$ & $\{\sink\}$ \\
 & $\delta(Hp \!\cdot\! K,K)$
   & $=$ & $\{\sink\}$ \quad\text{if}~ $HpK \in \sfl(\RbHK)$ \\
\hline
$H$-tree
 & $\delta(Hp,x)$
   & $=$ & $\left\{ \begin{array}{l}
            \{Hpx\} \quad\text{if}~ Hpx \in \ppl(\RbH) \\
            \{\sink\} \quad\text{if}~ Hpx \in \sfl(\RbH)
            \end{array} \right.$ \\
\hline
$K$-tree
 & $\delta(xqK,x)$
   & $=$ & $\{qK\} \quad\text{if}~ xqK \in \psl(\RbK)$ \\
\hline
$HK$-tree
 & $\delta(Hp \!\cdot\! K,x)$
   & $=$ & $\{Hpx \!\cdot\! K\} \quad\text{if}~ Hpx \in \ppl(\RbHK)$ \\
\hline
in $\ubNG$
 & $\delta(p,x)$
   & $=$ & $\left\{ \begin{array}{l}
            \{\sink\} \quad\text{if}~ px=ul
                      \text{ for some } l \in \sfl(\RbG),~ \text{else} \\
            \{p_i \in \ppl(\RbG) ~|~ px = u_ip_i\}
            \cup \{xqK ~|~ xqK \in \sfl(\RbK)\}
            \end{array} \right.$ \\
\hline
otherwise
 & $\delta(s,a)$
   & $=$ & $\emptyset$ \\
\hline
\end{tabular}

\medskip
Table 3 ~:~ Transition function $\delta$ for double cosets automaton.
\end{center}
\end{table}

\medskip
The extended state transition function $\delta^* : S \times \Sigma^* \to 2^S$ 
is such that, for $t \in \bT$, the intersection of $\delta^*(\init,t)$ 
with $S \setminus \{\norm\}$ is non-empty if and only if 
$t$ is a word in $\Sigma^*$ which is not in $\bT$ or is reducible.

Just as we converted $\ubNG$ to $\ubAG$, 
we make $\ubN$ deterministic; take its complement; and minimize.
The language $\bL$ recognised by the resulting automaton $\ubA$
is ${\bT_+} - ({\bT_+} - \irrRCT) = \irrRCT$.
Hence, by Kleene's Theorem,  $\bL$ is regular.
Since $R^C$ is a complete rewriting system on $\bT$,
there exists a unique irreducible word in each class of $\bT$
with respect to $\stackrel{*}{\leftrightarrow}_R$.
Therefore the set $\irrRCT$ is bijective with
$\bT/ \stackrel{*}{\leftrightarrow}_{R} = \bL$.

The automaton $\ubA$ gives rise to a
system of right linear language equations with a unique solution,
which is a regular expression for the language $\bL$ accepted by the automaton.
The regular expression can be obtained by applying Arden's Theorem
to solve the language equations.
\end{proof}

Thus an automaton $\ubA$ can be constructed from the complete
double coset rewriting system and a regular expression
for the set of double cosets $\bL$ is obtained from solving
the language equations of the determinized and minimized complement of $\ubA$.
Section \ref{sec-examples} includes some examples of these automata.


\section{Logged Double Coset Rewriting} \label{sec-logged}

It is often useful to label the original rewrite rules
and record how they are used during Knuth-Bendix completion,
and then during the reduction of words.
One instance is when we require precise proof of a
particular equivalence in terms of the original data.

Suppose that $\alpha_{s}:s \to t$ is a rewrite rule.
Then we know that the rewrite $u s v \to u t v$ can be performed,
and a reasonable label for this is $u \alpha_{s} v$.
Similarly if $\alpha_{t}:t \to q$ then we may perform $\alpha_{s}$
followed by $\alpha_{t}$, rewriting $s \to q$,
which we choose to label $\alpha_{s} \bullet \alpha_{t}$.
Thus any sequence of rewrites may be recorded by a
combination of labels of the form:
$$
u_1 \alpha_{s_1} v_1 \bullet u_2 \alpha_{s_2} v_2 \bullet \cdots
\bullet u_n \alpha_{s_n} v_n.
$$
The algebra of recorded rewrites is formalised by observing the
sesquicategorical or 2-categorical structure.
Briefly, a 2-category consists of 0-cells (objects $\bullet_i$),
1-cells (arrows between objects $(w_{ij} : \bullet_i \to \bullet_j)$) and
2-cells (arrows between arrows $(\alpha : w_{ij} \Rightarrow w'_{ij})$),
with a category structure on the 1-cells
(arrow composition) and two compatible (by the interchange law)
category structures (horizontal and vertical composition) on the 2-cells,
which preserve sources and targets~\citep{maclane}.

In particular {\em logged double coset rewriting},
is formalised in terms of a 2-category whose 0-cells (vertices) and
1-cells (paths along arrows) are illustrated in the following graph.
$$
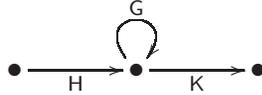
\xymatrix@C=3pc{
\bullet \ar[r]_{\;\bH\;}
   & \bullet \ar@(ul,ur)^{\bG}
             \ar[r]_{\;\bK\;}
      & \bullet \\}
$$
The generating 2-cells are
\begin{eqnarray*}
\alpha_{h_i} ~:~ H h_i \to H,
   && \text{ for each } h_i \in \XbH, \quad \mbox{ } \\
\alpha_{l_i} ~:~ \;l_i \to r_i, \quad
   && \text{ for each } (l_i,r_i) \in \RbG, \\
\text{and } \quad
\alpha_{k_i} ~:~ k_i K \to K,
   && \text{ for each } k_i \in \XbK. \, \quad \mbox{ }
\end{eqnarray*}

The vertical composition of 2-cells, written $\alpha \bullet \gamma$,
is the composition of the two rewrites 
when the source of the second coincides with the target of the first. 
We may also ``whisker'' the 2-cells with suitable 1-cells.
For example, the 2-cell $\alpha_{l_i}$ may be whiskered by
$Hu$ on the left and $v$ on the right to obtain a 2-cell
$(Hu \alpha_{l_i} v : Hul_iv \to Hur_iv)$.
Finally, we require the interchange law,              
so that it does not matter in which order we combine
a pair of 2-cells which rewrite non-overlapping parts of a string.
This defines the horizontal composition,
$$
\alpha \circ \beta
~:=~
\alpha\ \src(\beta) \bullet \tgt(\alpha) \beta
~=~
\src(\alpha) \beta \bullet \alpha\ \tgt(\beta), 
$$
and corresponds to the fact that if rewrite rules do not overlap
on a string then it does not matter which one we apply first.
Note that whiskering is equivalent to composing with identity 2-cells, 
$u \alpha v = 1_u \circ \alpha \circ 1_v$ 
where $(1_u : u \Rightarrow u)$ for $u \in \Sigma$, 
and that $\alpha v \circ \beta = \alpha \circ v\beta$. 

Logged rewriting for monoid presentations is explained in detail
in~\citep{paper11} and we shall only recall the key ideas here.

The essential difference between the logged version
and the standard Knuth-Bendix algorithm is the level of detail it records.
If we have an overlap which introduces a new rule,
we require an expression for the new rule in terms of the original labels.
Often we will be in the situation where $w$ reduces to $w_1$
by one sequence $\beta_1$ of 2-cells
and to $w_2$ by another sequence $\beta_2$.
Assuming $w_2$ is the larger string, we add in the rule $w_2 \to w_1$,
and note that this relation can be achieved by ``un-reducing'' $w_2$ to $w$
and then reducing $w$ to $w_1$.
This ``un-reducing'' is more formally known as an
\emph{inverse derivation}
and gives the vertical composition of 2-cells a groupoid structure.
In this situation we add the 2-cell $\beta_2^{-1} \bullet \beta_1$ 
at the same time as the rule $w_2 \to w_1$.

Assume that, using these methods,
we obtain a complete, logged, rewrite system for the double cosets,
which means that we have a 2-cell associated with each of the rules.
Suppose now that we have two group elements $g_1$ and $g_2$,
represented by strings $w_1$ and $w_2$ in the free monoid,
so that $\theta(w_1)=g_1$ and $\theta(w_2)=g_2$.
We can determine whether or not $g_1$ and $g_2$ lie within the same
double coset by rewriting $Hw_1K$ and $Hw_2K$.
If they both reduce to the same string $HzK$,
then we can examine the logs of the reductions to find
$h_1, \ldots, h_m \in \XbH$ and $k_1, \ldots, k_n \in \XbK$ such that
\quad
\begin{equation} \label{eq-logs}
\theta(h_1^{\ep_1}) \cdots \theta(h_n^{\ep_n})\; g_1\;
\theta(k_1^{\ep_{n+1}}) \cdots \theta(k_m^{\ep_{n+m}}) ~=~ g_2.
\end{equation}
The following section includes examples of this computation.

In the case of left or right cosets, 
the logs of the complete rewriting system may be used to derive 
a presentation for the subgroup itself.
In other situations the logs and particularly the logs of 
circular rewrites (endorewrites) have more interesting interpretations 
and applications~\citep{idrel}.
In the double coset case we can make the following observations.
We use $E(w)$ to denote the set of endorewrites of the string $w$,
the set of 2-cells associated to rewrites of $w$ back to itself.

\begin{enumerate}[i)]
\item
The sets $E(Hw)$ are bijective with each other for all $w$ in $X^*$.
These give information about the group $\bH$ in the form
of a presentation~\citep{RS-alt}.
\item
Similarly, the sets $E(wK)$ give a presentation of $\bK$.
\item
The sets $E(w)$ are all bijective, and these give generators for the
module of identities among relations for the group 
$\bG$~\citep{paper11,idrel}.
\item
The sets $E(HwK)$ are not all bijective in general.
However, in the case that $\theta(w_1)=\theta(w_2)$,
there is a bijection between the sets $E(Hw_1K)$ and $E(Hw_2K)$.
In general each endorewrite of this type gives us
information regarding the relationship between the subgroups
$\bH$ and $\bK$ within $\bG$.
Generators for the subgroup $\bH \cap \bK$ 
can certainly be obtained in this way.
\end{enumerate}


\section{Examples} \label{sec-examples}

The examples given below were calculated using a prototype package
{\sf kan}, available from~\citep{kan}.
This is a collection of {\sf GAP} functions which are
designed to tackle rewriting problems by translating them to a
categorical framework;
using a generalised Knuth-Bendix type algorithm to solve the problem;
and then translating back into the format
appropriate for the structure in question.
The double coset functions in {\sf kan} make use of functions from
the {\sf GAP} package {\sf automata}~\citep{automata}.
The main purpose of these examples is to demonstrate the methods we have
presented and the fact that they can be applied to a wider class of
problems than could previously be computed, rather than to be technically
impressive.

\begin{exmp}[A finite double coset rewriting system] \mbox{ }\\
\emph{Let $\bG$ be the free group on generators $\{a,b\}$
and let $\bH = \langle a^6 \rangle,~ \bK = \langle a^4 \rangle$.
(Varying the powers of $a$ gives a family of examples of this type.)
The double coset $\bH\bK$ contains $a^2 = a^{\gcd(6,4)}$,
and it is clear that one set of double coset representatives is
{\small $$
\{~ \bH \bK,~ \bH a \bK,~ \bH a^i b a^j \bK,~ \bH a^i bub a^j \bK   
~\mid~
0 \leqslant i \leqslant 5,~ 0 \leqslant j \leqslant 3,~
u \in \{a,b\}^* ~\}.
$$
}The initial set of rules is
{\small $$
\{~ (Aa,\id),~ (aA,\id),~ (Bb,\id),~ (bB,\id), (Ha^6,H),~ (a^4K,K) ~\}.
$$
}After completion, the last two rules are replaced by
{\small $$
\{ (Ha^4,HA^2), (HA^3,Ha^3), (a^3K,AK), (A^2K,a^2K),
(Ha^2K,HK), (HAK,HaK) \}.
$$
}Note that two $HK$-rules appear, reflecting the fact that
$\bH \vee \bK = \langle a^2 \rangle$.
The non-deterministic automaton $\ubN$ has $22$ states,
\begin{center}
\begin{tabular}{l}
{\small
$\{\init,\norm,\sink\}~\cup~\{\id,a,A,b,B\}~\cup~\{H,Ha,Ha^2,Ha^3,HA,HA^2\}$}
\\
\qquad
{\small
$\cup~\{K,aK,a^2K,AK\}~\cup~
\{H \!\cdot\! K,Ha \!\cdot\! K,Ha^2 \!\cdot\! K,HA \!\cdot\! K\}.$}
\end{tabular}
\end{center}
Determinizing $\ubN$ gives an automaton with $24$ states which,
after complementation and minimization, reduces to a deterministic
automaton with $15$ states and transitions shown in Table 4 
(where $\mathbb{N},\mathbb{S},\mathbb{I}$ correspond to $\norm,\sink,\init$). 
}

\medskip
\begin{table}[ht]  \label{table-ex7}
\begin{center}
{\small
\begin{tabular}{|c|ccccccccccccccc|}
\hline
    & $1$ & $2$ & $3$ & $4$ & $5$ & $6$ & $7$ & $8$ & $9$ & $10$
    & $11$ & $12$ & $13$ & $14$ & $15$ \\
\hline
    & $\mathbb{N}$ & $\mathbb{S}$ & $B$ & $\mathbb{I}$ & $a$ & $\id$ & $a^3$
    & $b$ & $BA$ & $a^2$ & $ba^3$ & $ba^2$ & $ba$ & $A$ & $BA^2$ \\
\hline
$H$ & $2$ & $2$ & $2$ & $6$ & $2$ & $2$ & $2$ & $2$ & $2$ & $2$
    & $2$ & $2$ & $2$ & $2$ & $2$ \\
$K$ & $2$ & $2$ & $1$ & $2$ & $1$ & $1$ & $2$ & $1$ & $1$ & $2$
    & $2$ & $1$ & $1$ & $2$ & $2$ \\
$a$ & $2$ & $2$ & $13$ & $2$ & $10$ & $5$ & $2$ & $13$ & $2$ & $7$
    & $11$ & $11$ & $12$ & $2$ & $2$ \\
$A$ & $2$ & $2$ & $9$ & $2$ & $2$ & $14$ & $2$ & $9$ & $15$ & $2$
    & $2$ & $2$ & $2$ & $7$ & $15$ \\
$b$ & $2$ & $2$ & $2$ & $2$ & $8$ & $8$ & $8$ & $8$ & $8$
    & $8$ & $8$ & $8$ & $8$ & $8$ & $8$ \\
$B$ & $2$ & $2$ & $3$ & $2$ & $3$ & $3$ & $3$ & $2$
    & $3$ & $3$ & $3$ & $3$ & $3$ & $3$ & $3$ \\
\hline
\end{tabular}}

\medskip
Table 4 ~:~ Minimal double coset automaton for Example 7.
\end{center}
\end{table}
\end{exmp}

\bigskip
\begin{exmp}[The Trefoil Group]\mbox{ }\\
\emph{This is an example in which the group has 
a finite rewriting system,
but the double coset system is infinite.
Starting with an initial monoid presentation 
with rules} 
$$
{\small 
[~\alpha_1 = (Xx, \id),\ 
\alpha_2 = (xX, \id),\ 
\alpha_3 = (Yy, \id),\ 
\alpha_4 = (yY, \id),\ 
\alpha_5 = (x^3, y^2)~], 
}$$
\emph{the fundamental group $\calT = \langle x,y ~|~ x^3=y^2 \rangle$
of the trefoil knot has a complete rewriting system with six logged rules 
shown in Table 5.} 

\medskip
\begin{table}[ht]  \label{table-ex8}
\begin{center}
{\small
\begin{tabular}{|c|l|}
\hline
rule & \qquad label\\
\hline
$(Yy, \id)$ & $\alpha_3$ \\
$(yY, \id)$ & $\alpha_4$ \\
$(x^3, y^2)$ & $\alpha_5$ \\
$(y^2x, xy^2)$ 
  & $\alpha_6 = (\alpha_5^{-1}x)\bullet(x\alpha_5)$ \\
$(X, x^2Y^2)$ 
  & $\alpha_7 = (X\alpha_4^{-1})\bullet(Xy\alpha_4^{-1}Y)
     \bullet(X\alpha5^{-1}Y^2)\bullet(\alpha_1x^2Y^2)$ \\
$(Yx, yxY^2)$ 
  & $\alpha_8 = (Yx\alpha_4^{-1})\bullet(Yxy\alpha_4^{-1}Y)
    \bullet(Yx\alpha_5^{-1}Y^2)\bullet(Y\alpha_5xY^2)\bullet(\alpha_3yxY^2)$ \\
\hline
\end{tabular}}

\medskip
Table 5 ~:~ Logged rewrite rules for the trefoil monoid.
\end{center}
\end{table}

\medskip
\emph{The ordering used here is the wreath product order with $X>x>Y>y$.
A group version of these logged rules is given in~\citep{idrel}.}

\emph{The non-deterministic automaton
$\ubN_{\calT}$ has $7$ states,
and there are $12$ states in the determinized automaton,
reducing to $7$ states on minimization.
The automaton $\ubAT$ is pictured in Figure 2. 
(For clarity, the sink state and transitions to it have been excluded.
All states are accepting, so double circles have been omitted.)}

\medskip
\begin{figure}[htb] \label{fig-trefoil-dfa}
\begin{center}
\input{jsc-3474-fig2.pstex_t}

\medskip
Figure 2 ~:~ Sketch of the word acceptor  $\ubAT$ for $\calT$
\end{center}

\medskip
\end{figure}
\emph{A regular expression for the language
accepted by $\ubAT$ is}
$$
(1+y)x(yx+xyx)^*(1+x)(y^*+Y^+) + (y^*+Y^+).
$$

\emph{We consider subgroups $\bH = \langle x \rangle$ and
$\bK = \langle y \rangle$.
This is an example which apparently cannot be computed using
algorithms previously available.
The double coset rewriting system initially requires the additional rules}
{\small $$
[~\beta_1 = (HX, H),\ 
\beta_2 = (Hx, H),\ 
\beta_3 = (YK, K),\ 
\beta_4 = (yK, K)~]. 
$$
}\emph{The {\sf kan} package includes a \emph{limited} 
version of the Knuth-Bendix functions which stop the calculation 
after a specified number $\lambda$ of rules have been added to the system.
Subsets of rules which involve either $H$, or $K$, or both 
may then be extracted.
Adopting the wreath product order with $K>H>X>x>Y>y$ we find that
$\beta_1$ reduces, and there are no additional $K$-rules or $HK$-rules.
As $\lambda$ increases we obtain an increasing number of rules from 
the following infinite set:
{\small $$
\{~ (Hwy^2,Hw),~ (HwY, Hwy) \quad|\quad
w ~\text{ is any word in }~ \{yx,yx^2\}^* ~\}.
$$
}Note that the words $yx$ and $yx^2$ label directed cycles in Figure 2. 
It is straightforward to verify that, if we add all these rules, 
the system is complete, despite the fact that it is infinite:
the $H$ acting as a tag on the left restricts the possible overlaps
severely and only a few cases need be checked.} 

\emph{Also as $\lambda$ increases, the left-hand sides of the rules may be 
used to form the three trees in the double coset automaton of Figure 1.
Applying the construction of Theorem \ref{thm-dcwa}, 
we obtain a sequence of minimized automata $\{\ubAl\}$ 
which rapidly converges to the automaton depicted in Figure 3. 
Indeed, with $\lambda=10$, we obtain a sufficient set of 
three $H$-rules and two $K$-rules, namely  
$\{\beta_2,\beta_3,\beta_4\}$ together with 
{\small
\begin{center}
\begin{tabular}{|c|l|}
\hline
rule & \qquad label\\
\hline
$\beta_5 = (Hy^2, H)$ 
  & $H\alpha_5^{-1} \bullet \beta_2x^2 \bullet \beta_2x 
                \bullet \beta_2$ \\
$\beta_6 = (HY, Hy)$ 
  & $\beta_5^{-1}Y \bullet Hy\alpha_4$ \\
\hline
\end{tabular}
\end{center}
}}

\medskip
\begin{figure}[hbt] \label{fig-trefoil-dc}
\begin{center}
\input{jsc-3474-fig3.pstex_t}

\medskip
Figure 3 ~:~ Double coset automaton for the trefoil group. 
\end{center}

\medskip
\end{figure}

\emph{The normal forms can be read straight off
the automaton:
{\small $$
\{~ HwK  \quad|\quad
w ~\text{ is any word in }~ \{yx,yx^2\}^* ~\}.
$$
}}\emph{For an example of logged reduction, consider the double coset
$\bH Y \bK$, where $Y$ is a normal form for $\ubAT$.
Applying $H\beta_3$ we obtain immediately that $HYK \to HK$.
Alternatively, we may apply 
\small{$$
\beta_6 K \bullet H \beta_4 
   ~=~ \beta_2^{-1}YK \bullet \beta_2^{-1}xYK \bullet 
       \beta_2^{-1}x^2YK \bullet H\alpha_5K \bullet Hy\alpha_4K 
       \bullet H \beta_4, 
$$}which gives successive rewrites
\small{$$
HYK \to HxYK \to Hx^2YK \to Hx^3YK \to Hy^2YK \to HyK \to HK.
$$}
Applying equation (\ref{eq-logs}) to these two reductions 
(where $w_1 = Y$ and $w_2 = \id$) gives} 
\begin{eqnarray*}
\theta(w_1)\theta(Y)^{-1}
  & ~=~ &  y^{-1}y ~=~ 1_{\bG} ~=~ \theta(w_2), \\
(\theta(x))^3\theta(w_1)\theta(y)^{-1}
  & ~=~ &  \theta(x^3Y^2) ~=~ 1_{\bG} ~=~ \theta(w_2), \\
\end{eqnarray*}
\emph{and so we have obtained an endorewrite 
$H \beta_3^{-1} \bullet \beta_6 K \bullet H \beta_4 ~\in~ E(HK)$.}
\end{exmp}

\bigskip
\begin{exmp}[A group with an infinite rewriting system]\mbox{ }\\
\emph{When the group $\bG$ has an infinite rewriting system
the double coset rewriting system will also be infinite.
In this case it may be possible to use the package {\sf KBMAG}
to compute a word acceptor for $\bG$.
In the {\sf kan} package the finite state automaton provided by {\sf KBMAG} 
is converted to a deterministic automaton in the format used by the 
{\sf automata} package, and then included as the $\ubNG$ part of 
the double coset automaton shown in Figure 1.
It is still necessary to find appropriate sets of rules 
$\RbH, \RbK$ and $\RbHK$ and, since $\RbG$ is infinite, 
the limited Knuth-Bendix functions should again be used.
An interactive use of the package is required:
experimenting with different limits gives partial results
from which we may be able to deduce exact answers.} 

\emph{In the following example we take $\bG$ to be the group with generators
$\{a,b\}$ and relators $[a^3,b^3,(ab)^3]$.
The normal forms of the monoid elements are strings alternating in
$a$ or $A$ with $b$ or $B$.
Not all such strings are irreducible, for example 
$bab \to ABA$ and $abaB \to BAb$.
The automatic structure computed by {\sf KBMAG} has a word acceptor
with $17$ states.}

\emph{If we take $\bH$ to be generated by $[ab]$ 
and $\bK$ to be generated by $[ba]$, 
we find that all three additional sets of rules are infinite: 
{\small 
\begin{eqnarray*}
 \RbH &~=~& \{~ (Hab,H),~ (HaB,Hb),~ (H(bA)^nB,H(Ab)^na), n \geqslant 0 ~\}, \\
 \RbK &~=~& \{~ (baK,K),~ (BaK,bK),~ (B(Ab)^nK,a(bA)^nK), n \geqslant 0 ~\}, \\
\RbHK &~=~& \{~ (Hb(Ab)^nK,H(Ab)^nAK), n \geqslant 0 ~\}. 
\end{eqnarray*}
}The sequence $\{\ubAl\}$ of minimized automata exhibits an increasing 
number of states, reaching $48$ at $\lambda=200$.  
On inspection, we find that these automata have a common set of states 
(the right-hand half in Figure 4),
and two chains which gradually increase in length.
Taking $\ubA = \lim_{\lambda\to\infty} \ubAl$, 
these chains shrink to two-state loops, and
we obtain the $19$ state automaton shown in Figure 4, 
where the $\norm$ and $\sink$ states have been omitted, 
and states shown with a double circle are those having 
a $K$-transition to $\norm$.} 

\medskip
\begin{figure}[htb] \label{fig-ex3}
\begin{center}
\input{jsc-3474-fig4.pstex_t}

\medskip
Figure 4 ~:~ Automaton $\ubA$ for Example 9. 
\end{center}

\medskip
\end{figure}
\emph{The language recognised by $\ubA$ is}
{\small
$$
H \left(~ 
a + (\id + A)(bA)^* + AB(aB)^*A + b(aB)^+A(bA)^* + A(bA)^*(Ba)^*b 
~\right) K.
$$}  
\end{exmp}

\section{Induced actions and left Kan extensions}
\label{sec-kan}

The algorithms we have presented arose as part of a programme 
of applying categorical constructions in computer algebra 
so as to allow increased flexibility and a wider range of applications.
It is worth recalling~\citep{maclane} in noting that our goal
\emph{is not the reduction of the familiar to the unfamiliar, but the
extension of the familiar to cover many more cases}. 
In this section, we shall demonstrate this by showing how the double coset
problem is an instance of the much more general construction of 
\emph{induced actions of monoids and categories}.  
Our aim also is to advertise this construction, 
which has many uses apart from those given here 
(see, for example, \citep{paper2}).

Let $F: \bM \to \bN$ be a morphism of monoids, 
and let $\bM,\bN$ act on sets $X,Y$ respectively. 
These are right actions, and we use the notation $x^m,y^n$.
An \emph{$F$-morphism} $\ep: X \to Y$ is a function of sets such that
$\ep(x^m) = (\ep x)^{F(m)}$ for all $m \in \bM, x \in X$. 
In the case $F = 1_{\bN}$ we call $\ep$ an $\bN$-morphism. 
The $F$-morphism $\ep$ is said to be \emph{universal} if, 
for any action of $\bN$ on a set $Z$ and any $F$-morphism $\phi: X \to Z$, 
there is a unique $\bN$-morphism $\psi: Y \to Z$ 
such that $\psi \circ \ep = \phi$. 
When $\ep$ is universal we say that the action of $\bN$ on $Y$ is 
\emph{induced} from $X$ by $F$, and we write $Y=F_*(X)$.

It is easy to construct this $Y$ from the action $X \times \bM \to X$
and the morphism $F$. 
We let $Y' = X \times \bN$ with $\bN$-action $(x,n)^{n'} := (x,nn')$, 
and define $Y$ to be $Y'$ factored by the
equivalence relation $\approx$ generated by
\begin{equation} \label{eq-actionequiv}
(x^m,n) ~\approx~ (x,F(m)n), \quad n \in \bN,~ m \in \bM,~ x \in X, 
\end{equation}
with $\ep : X \to Y$ mapping $x$ to the equivalence class of $(x,1_{\bN})$.

A common example is where $\bM$ is a subgroup of a group $\bN$, 
$F$ is the inclusion morphism, and $X$ is a singleton. 
Then $F_*(X)$ can easily be identified with $\bM\backslash\bN$, 
the set of left cosets of $\bM$ in $\bN$, and the usual $\bN$-action.
Note also that we immediately have an extension of the problem 
of determining cosets to that of extending an action of the subgroup $\bM$ 
on a set $X$ to an action of $\bN$ on a new set $F_*(X)$.

This notion of induced action is easily extended to the case where
$\bM,\bN$ are small categories, $F : \bM\to\bN$ is a functor, 
and the action of $\bM$ is given by a functor 
$X : \bM \to \sets$ (see below). 
In~\citep{heyworth-thes,paper2} string rewriting for a monoid
presentation of a monoid $\bN$ was generalised to string rewriting for
an induced action of categories, given a presentation of the data for this. 
Our string rewriting procedure for double cosets is a special case 
of this general form of string rewriting for induced actions of categories.

Note that induced actions are also well known in category theory
as left Kan extensions, and have many applications under that
terminology. In that setting, it is often convenient to describe a
choice of left Kan extensions for all actions as giving a functor
of functor categories
$$
F_* ~:~ (\sets)^{\bM} \to (\sets)^{\bN}
$$ 
which is left adjoint to the standard functor
$$
F^* ~:~ (\sets)^{\bN} \to (\sets)^{\bM}
$$ 
given by composition with $F$.
An implication of this is that $F_*$ preserves colimits of
actions, but we do not pursue this theme.

The construction of the functor $Y = F_*(X) : \bN \to \sets$  
in the category case is an easy generalisation of that for monoids, 
but with appropriate attention to the objects of the categories concerned. 
For $b \in \ob(\bN)$ let $Y'(b)$ be the disjoint union of the sets 
$X(c) \times \bN(F(c),b)$ for all $c \in \ob(\bM)$. 
On $Y'(b)$ we impose the equivalence relation generated by 
equation (\ref{eq-actionequiv}) to give $Y(b)$ as the quotient. 
Note that now $x \in X(c),~ x^m \in X(c'),~ m \in \bM(c,c')$ 
and $n \in \bN(c',b)$, as in Figure 5.
The action of $\bN$ is induced by $(x,n)^{n'} = (x,nn')$ as before.
This construction is known in category theory as that of a coend.

\medskip
\begin{figure}[htb] \label{fig-gpd-act}
\begin{center}
\input{jsc-3474-fig5.pstex_t}

\medskip
Figure 5 ~:~ Action induced from $X$ by $F$. 
\end{center}

\medskip
\end{figure}
We now apply this construction to the double coset problem. 
Let $\Gamma$ be a set with commuting right $\bH$- and $\bK$-actions, 
so that $(\gamma^h)^k = (\gamma^k)^h$ for all
$\gamma\in\Gamma, h\in\bH, k\in\bK$.
Intuitively we prefer to think of the $\bH$-action as a left action, 
defining ${}^{h^{-1}}\gamma := \gamma^h$.
We define a category $\bM = \bH\circto\bK$ with objects $\{1,2\}$ 
and the following elements:
$$
\bM(1,1) = \bH,\quad 
\bM(2,2) = \bK,\quad 
\bM(1,2) = \bH\times\bK,\quad 
\bM(2,1) = \emptyset.
$$
Composition in $\bM$ is given by the usual multiplication in $\bH,\bK$ 
and by
$$
h_2(h_1,k_1)k_2 ~=~ (h_2h_1,k_1k_2),
\qquad\text{so that}\qquad
(h,k) ~=~ h(1,1)k. 
$$
For the $\bM$-action $X$ we take $X(1) = \Gamma\times\{1\}$ and 
$X(2) = \Gamma\times\{2\}$ to be copies of $\Gamma$.  
For the morphisms we define 
$$
X(h)(\gamma,1) = ({}^{h^{-1}}\gamma,1),\quad 
X(k)(\gamma,2) = (\gamma^k,2),\quad 
X(h,k)(\gamma,1) = ({}^{h^{-1}}\gamma^k,2), 
$$
so $X(1,1)$ is the isomorphism between copies of $\Gamma$  
mapping $(\gamma,1)$ to $(\gamma,2)$. 

\begin{prop} 
If $\bN$ is the trivial category with one object $0$ 
and $1_0$ the only morphism, 
$X$ is the $\bM$-action given above, 
and $F$ is the unique functor $\bM \to \bN$, 
then $F_*(X)$ may be identified with the set of orbits $\HGammaK$. 
\end{prop}
\begin{proof}
An $\bN$-action $Y$ is just a set $Y(0)$ and the identity
function $Y(1_0)$ on $Y(0)$.
The construction above gives 
$Y'(0) = (\Gamma\times\{1\}\times\{1_{0}\}) 
\sqcup (\Gamma\times\{2\}\times\{1_{0}\})$.
Applying the equivalence rule (\ref{eq-actionequiv}) with 
$m = h \in \bH,~ m = (1,1)$ and $m = k \in \bK$, we obtain 
$$
({}^{h^{-1}}\gamma,1,1_{0}) \approx (\gamma,1,1_{0}),~ 
(\gamma,2,1_{0}) \approx (\gamma,1,1_{0}),~ 
(\gamma^k,2,1_{0}) \approx (\gamma,2,1_{0}), 
$$
identifying the two copies of $\Gamma$ 
and modelling the $\bH,\bK$ actions on $\Gamma$. 
\end{proof}

When $\Gamma = \bG$ and the actions are $g^h:=h^{-1}g,~ g^k:=gk$, 
then $F_*(X)$ may be identified with the set of double cosets $\HGK$. 

To summarise, as with cosets and other problems in computational
algebra, double cosets are an instance of the general problem of
computing Kan extensions or induced actions of categories. 
By developing string rewriting for computing such Kan extensions, 
we therefore have a generic algorithm applicable to all of these problems.

\section{Conclusions}  \label{sec-comcon}

One of the nice outcomes of our results is that existing powerful
string rewriting software can immediately be applied to double
coset problems, provided that we have a presentation for the group
and generating sets for the subgroups. Of course the algorithm may
not terminate, but {it will} in all cases where there are
a finite number of cosets $\bH  g$ and $g \bK$, 
(where the existing enumerative methods could be used), 
and also in some cases where there are an infinite number of double cosets 
and enumerative methods are likely to fail.

We have developed procedures for logged string rewriting, which
via a 2-categorical structure, records all computations made from
the original presentations. This enables us { to compute
not only } whether two elements of a group lie within the same
double coset but also, {using logged rewriting,  to
produce a proof of this} in the case that they do. We expect to
release the {\sf GAP} functions to do both determining and proving
via logged string rewriting as a share package.

In this paper we apply our results only to group theory, but
as we showed above, they hold much more generally.
It might be interesting to see whether there are analogues to
the double coset problem in other structures such as monoids or algebras --
structures where we also already have Kan extension rewriting methods
available to us.


\bibliographystyle{elsart-harv}
\bibliography{jsc-3474}

\end{document}

%% file: jsc-3474-fig1.pstex_t
\begin{picture}(0,0)%
\includegraphics{jsc-3474-fig1.pstex}%
\end{picture}%
\setlength{\unitlength}{1776sp}%
\begingroup\makeatletter\ifx\SetFigFont\undefined%
\gdef\SetFigFont#1#2#3#4#5{%
  \reset@font\fontsize{#1}{#2pt}%
  \fontfamily{#3}\fontseries{#4}\fontshape{#5}%
  \selectfont}%
\fi\endgroup%
\begin{picture}(13974,6119)(-2861,-4897)
\put(3376,-811){\makebox(0,0)[lb]{\smash{{\SetFigFont{9}{10.8}{\familydefault}{\mddefault}{\updefault}{\color[rgb]{0,0,0}$x$}%
}}}}
\put(8776,-1411){\makebox(0,0)[lb]{\smash{{\SetFigFont{9}{10.8}{\familydefault}{\mddefault}{\updefault}{\color[rgb]{0,0,0}$K$-tree}%
}}}}
\put(901,-2911){\makebox(0,0)[lb]{\smash{{\SetFigFont{9}{10.8}{\familydefault}{\mddefault}{\updefault}{\color[rgb]{0,0,0}$HK$-tree}%
}}}}
\put(-1724,-3286){\makebox(0,0)[lb]{\smash{{\SetFigFont{9}{10.8}{\familydefault}{\mddefault}{\updefault}{\color[rgb]{0,0,0}$H$-tree}%
}}}}
\put(1276,764){\makebox(0,0)[lb]{\smash{{\SetFigFont{9}{10.8}{\familydefault}{\mddefault}{\updefault}{\color[rgb]{0,0,0}$\norm$}%
}}}}
\put(-1499,-4036){\makebox(0,0)[lb]{\smash{{\SetFigFont{9}{10.8}{\familydefault}{\mddefault}{\updefault}{\color[rgb]{0,0,0}$Hp$}%
}}}}
\put(-1499,-2161){\makebox(0,0)[lb]{\smash{{\SetFigFont{9}{10.8}{\familydefault}{\mddefault}{\updefault}{\color[rgb]{0,0,0}$H\id$}%
}}}}
\put(7501,-2986){\makebox(0,0)[lb]{\smash{{\SetFigFont{9}{10.8}{\familydefault}{\mddefault}{\updefault}{\color[rgb]{0,0,0}$K$}%
}}}}
\put(8950,-211){\makebox(0,0)[lb]{\smash{{\SetFigFont{9}{10.8}{\familydefault}{\mddefault}{\updefault}{\color[rgb]{0,0,0}$xqK$}%
}}}}
\put(1126,-2236){\makebox(0,0)[lb]{\smash{{\SetFigFont{9}{10.8}{\familydefault}{\mddefault}{\updefault}{\color[rgb]{0,0,0}$H\!\cdot\!K$}%
}}}}
\put(3376,-211){\makebox(0,0)[lb]{\smash{{\SetFigFont{9}{10.8}{\familydefault}{\mddefault}{\updefault}{\color[rgb]{0,0,0}$\id$}%
}}}}
\put(9001,-2761){\makebox(0,0)[lb]{\smash{{\SetFigFont{9}{10.8}{\familydefault}{\mddefault}{\updefault}{\color[rgb]{0,0,0}$\id K$}%
}}}}
\put(5825,-2761){\makebox(0,0)[lb]{\smash{{\SetFigFont{9}{10.8}{\familydefault}{\mddefault}{\updefault}{\color[rgb]{0,0,0}$\sink$}%
}}}}
\put(-1949,-136){\makebox(0,0)[lb]{\smash{{\SetFigFont{9}{10.8}{\familydefault}{\mddefault}{\updefault}{\color[rgb]{0,0,0}$\init$}%
}}}}
\put(1051,-3436){\makebox(0,0)[lb]{\smash{{\SetFigFont{9}{10.8}{\familydefault}{\mddefault}{\updefault}{\color[rgb]{0,0,0}$Hp\!\cdot\!K$}%
}}}}
\put(1501,-4711){\makebox(0,0)[lb]{\smash{{\SetFigFont{9}{10.8}{\familydefault}{\mddefault}{\updefault}{\color[rgb]{0,0,0}$x$}%
}}}}
\put(3301,-4261){\makebox(0,0)[lb]{\smash{{\SetFigFont{9}{10.8}{\familydefault}{\mddefault}{\updefault}{\color[rgb]{0,0,0}$K$}%
}}}}
\put(3526,989){\makebox(0,0)[lb]{\smash{{\SetFigFont{9}{10.8}{\familydefault}{\mddefault}{\updefault}{\color[rgb]{0,0,0}$K$}%
}}}}
\put(-1799,-1111){\makebox(0,0)[lb]{\smash{{\SetFigFont{9}{10.8}{\familydefault}{\mddefault}{\updefault}{\color[rgb]{0,0,0}$H$}%
}}}}
\put(5200,-400){\makebox(0,0)[lb]{\smash{{\SetFigFont{9}{10.8}{\familydefault}{\mddefault}{\updefault}{\color[rgb]{0,0,0}$p$}%
}}}}
\put(5550,-1525){\makebox(0,0)[lb]{\smash{{\SetFigFont{9}{10.8}{\familydefault}{\mddefault}{\updefault}{\color[rgb]{0,0,0}$px$}%
}}}}
\put(7351,914){\makebox(0,0)[lb]{\smash{{\SetFigFont{9}{10.8}{\familydefault}{\mddefault}{\updefault}{\color[rgb]{0,0,0}$x$}%
}}}}
\put(451,-436){\makebox(0,0)[lb]{\smash{{\SetFigFont{9}{10.8}{\familydefault}{\mddefault}{\updefault}{\color[rgb]{0,0,0}$H$}%
}}}}
\put(-224,-1411){\makebox(0,0)[lb]{\smash{{\SetFigFont{9}{10.8}{\familydefault}{\mddefault}{\updefault}{\color[rgb]{0,0,0}$H$}%
}}}}
\put(3676,-1336){\makebox(0,0)[lb]{\smash{{\SetFigFont{9}{10.8}{\familydefault}{\mddefault}{\updefault}{\color[rgb]{0,0,0}$x$}%
}}}}
\put(3751,-2461){\makebox(0,0)[lb]{\smash{{\SetFigFont{9}{10.8}{\familydefault}{\mddefault}{\updefault}{\color[rgb]{0,0,0}$\ubNG$}%
}}}}

\put(4576,-1036){\makebox(0,0)[lb]{\smash{{\SetFigFont{9}{10.8}{\familydefault}{\mddefault}{\updefault}{\color[rgb]{0,0,0}$x$}%
}}}}
\put(5626,-961){\makebox(0,0)[lb]{\smash{{\SetFigFont{9}{10.8}{\familydefault}{\mddefault}{\updefault}{\color[rgb]{0,0,0}$x$}%
}}}}
\end{picture}%

%% file: jsc-3474-fig2.pstex_t
\begin{picture}(0,0)%
\includegraphics{jsc-3474-fig2.pstex}%
\end{picture}%
\setlength{\unitlength}{2368sp}%
\begingroup\makeatletter\ifx\SetFigFont\undefined%
\gdef\SetFigFont#1#2#3#4#5{%
  \reset@font\fontsize{#1}{#2pt}%
  \fontfamily{#3}\fontseries{#4}\fontshape{#5}%
  \selectfont}%
\fi\endgroup%
\begin{picture}(6996,3987)(-86,-2806)
\put(5776,-2650){\makebox(0,0)[lb]{\smash{{\SetFigFont{7}{8.4}{\rmdefault}{\mddefault}{\updefault}{\color[rgb]{0,0,0}$y$}%
}}}}
\put(1651,-700){\makebox(0,0)[lb]{\smash{{\SetFigFont{7}{8.4}{\rmdefault}{\mddefault}{\updefault}{\color[rgb]{0,0,0}$x$}%
}}}}
\put(4426,-700){\makebox(0,0)[lb]{\smash{{\SetFigFont{7}{8.4}{\rmdefault}{\mddefault}{\updefault}{\color[rgb]{0,0,0}$x$}%
}}}}
\put(3098,-878){\makebox(0,0)[lb]{\smash{{\SetFigFont{7}{8.4}{\rmdefault}{\mddefault}{\updefault}{\color[rgb]{0,0,0}$x$}%
}}}}
\put(2900,-1750){\makebox(0,0)[lb]{\smash{{\SetFigFont{7}{8.4}{\rmdefault}{\mddefault}{\updefault}{\color[rgb]{0,0,0}$x$}%
}}}}
\put(3301,-1750){\makebox(0,0)[lb]{\smash{{\SetFigFont{7}{8.4}{\rmdefault}{\mddefault}{\updefault}{\color[rgb]{0,0,0}$y$}%
}}}}
\put(3076,-2650){\makebox(0,0)[lb]{\smash{{\SetFigFont{7}{8.4}{\rmdefault}{\mddefault}{\updefault}{\color[rgb]{0,0,0}$y$}%
}}}}
\put(5750,-886){\makebox(0,0)[lb]{\smash{{\SetFigFont{7}{8.4}{\rmdefault}{\mddefault}{\updefault}{\color[rgb]{0,0,0}$x^2$}%
}}}}
\put(376,-886){\makebox(0,0)[lb]{\smash{{\SetFigFont{7}{8.4}{\rmdefault}{\mddefault}{\updefault}{\color[rgb]{0,0,0}$\id$}%
}}}}
\put(1651,164){\makebox(0,0)[lb]{\smash{{\SetFigFont{7}{8.4}{\rmdefault}{\mddefault}{\updefault}{\color[rgb]{0,0,0}$Y$}%
}}}}
\put(3226, 14){\makebox(0,0)[lb]{\smash{{\SetFigFont{7}{8.4}{\rmdefault}{\mddefault}{\updefault}{\color[rgb]{0,0,0}$Y$}%
}}}}
\put(4426,164){\makebox(0,0)[lb]{\smash{{\SetFigFont{7}{8.4}{\rmdefault}{\mddefault}{\updefault}{\color[rgb]{0,0,0}$Y$}%
}}}}
\put(3826,914){\makebox(0,0)[lb]{\smash{{\SetFigFont{7}{8.4}{\rmdefault}{\mddefault}{\updefault}{\color[rgb]{0,0,0}$Y$}%
}}}}
\put(1651,-1861){\makebox(0,0)[lb]{\smash{{\SetFigFont{7}{8.4}{\rmdefault}{\mddefault}{\updefault}{\color[rgb]{0,0,0}$y$}%
}}}}
\put(4426,-1936){\makebox(0,0)[lb]{\smash{{\SetFigFont{7}{8.4}{\rmdefault}{\mddefault}{\updefault}{\color[rgb]{0,0,0}$y$}%
}}}}
\put(3076,914){\makebox(0,0)[lb]{\smash{{\SetFigFont{7}{8.4}{\rmdefault}{\mddefault}{\updefault}{\color[rgb]{0,0,0}$Y$}%
}}}}
\put(4426,-2761){\makebox(0,0)[lb]{\smash{{\SetFigFont{7}{8.4}{\rmdefault}{\mddefault}{\updefault}{\color[rgb]{0,0,0}$y$}%
}}}}
\put(6526,-2650){\makebox(0,0)[lb]{\smash{{\SetFigFont{7}{8.4}{\rmdefault}{\mddefault}{\updefault}{\color[rgb]{0,0,0}$y^2$}%
}}}}
\end{picture}%

%% file: jsc-3474-fig3.pstex_t
\begin{picture}(0,0)%
\includegraphics{jsc-3474-fig3.pstex}%
\end{picture}%
\setlength{\unitlength}{1973sp}%
\begingroup\makeatletter\ifx\SetFigFont\undefined%
\gdef\SetFigFont#1#2#3#4#5{%
  \reset@font\fontsize{#1}{#2pt}%
  \fontfamily{#3}\fontseries{#4}\fontshape{#5}%
  \selectfont}%
\fi\endgroup%
\begin{picture}(6263,3858)(-2036,-2786)
\put(2026,800){\makebox(0,0)[lb]{\smash{{\SetFigFont{6}{7.2}{\rmdefault}{\mddefault}{\updefault}{\color[rgb]{0,0,0}$\norm$}%
}}}}
\put(376,-850){\makebox(0,0)[lb]{\smash{{\SetFigFont{6}{7.2}{\rmdefault}{\mddefault}{\updefault}{\color[rgb]{0,0,0}$\id$}%
}}}}
\put(-599,-700){\makebox(0,0)[lb]{\smash{{\SetFigFont{6}{7.2}{\rmdefault}{\mddefault}{\updefault}{\color[rgb]{0,0,0}$H$}%
}}}}
\put(2176,-700){\makebox(0,0)[lb]{\smash{{\SetFigFont{6}{7.2}{\rmdefault}{\mddefault}{\updefault}{\color[rgb]{0,0,0}$x$}%
}}}}
\put(950, 14){\makebox(0,0)[lb]{\smash{{\SetFigFont{6}{7.2}{\rmdefault}{\mddefault}{\updefault}{\color[rgb]{0,0,0}$K$}%
}}}}
\put(3350, 14){\makebox(0,0)[lb]{\smash{{\SetFigFont{6}{7.2}{\rmdefault}{\mddefault}{\updefault}{\color[rgb]{0,0,0}$K$}%
}}}}
\put(1051,-1850){\makebox(0,0)[lb]{\smash{{\SetFigFont{6}{7.2}{\rmdefault}{\mddefault}{\updefault}{\color[rgb]{0,0,0}$y$}%
}}}}
\put(3920,-850){\makebox(0,0)[lb]{\smash{{\SetFigFont{6}{7.2}{\rmdefault}{\mddefault}{\updefault}{\color[rgb]{0,0,0}$yx$}%
}}}}
\put(2176,-2650){\makebox(0,0)[lb]{\smash{{\SetFigFont{6}{7.2}{\rmdefault}{\mddefault}{\updefault}{\color[rgb]{0,0,0}$y$}%
}}}}
\put(3226,-2011){\makebox(0,0)[lb]{\smash{{\SetFigFont{6}{7.2}{\rmdefault}{\mddefault}{\updefault}{\color[rgb]{0,0,0}$y$}%
}}}}
\put(2851,-1561){\makebox(0,0)[lb]{\smash{{\SetFigFont{6}{7.2}{\rmdefault}{\mddefault}{\updefault}{\color[rgb]{0,0,0}$x$}%
}}}}
\put(-1499,-850){\makebox(0,0)[lb]{\smash{{\SetFigFont{6}{7.2}{\rmdefault}{\mddefault}{\updefault}{\color[rgb]{0,0,0}$\init$}%
}}}}
\end{picture}%

%% file: jsc-3474-fig4.pstex_t
\begin{picture}(0,0)%
\includegraphics{jsc-3474-fig4.pstex}%
\end{picture}%
\setlength{\unitlength}{1776sp}%
\begingroup\makeatletter\ifx\SetFigFont\undefined%
\gdef\SetFigFont#1#2#3#4#5{%
  \reset@font\fontsize{#1}{#2pt}%
  \fontfamily{#3}\fontseries{#4}\fontshape{#5}%
  \selectfont}%
\fi\endgroup%
\begin{picture}(13141,4167)(-6158,-4720)
\put(260,-1636){\makebox(0,0)[lb]{\smash{{\SetFigFont{5}{6.0}{\rmdefault}{\mddefault}{\updefault}{\color[rgb]{0,0,0}$b$}%
}}}}
\put(-524,-886){\makebox(0,0)[lb]{\smash{{\SetFigFont{5}{6.0}{\rmdefault}{\mddefault}{\updefault}{\color[rgb]{0,0,0}$1$}%
}}}}
\put(1276,-886){\makebox(0,0)[lb]{\smash{{\SetFigFont{5}{6.0}{\rmdefault}{\mddefault}{\updefault}{\color[rgb]{0,0,0}$2$}%
}}}}
\put(3076,-700){\makebox(0,0)[lb]{\smash{{\SetFigFont{5}{6.0}{\rmdefault}{\mddefault}{\updefault}{\color[rgb]{0,0,0}$a$}%
}}}}
\put(226,-700){\makebox(0,0)[lb]{\smash{{\SetFigFont{5}{6.0}{\rmdefault}{\mddefault}{\updefault}{\color[rgb]{0,0,0}$H$}%
}}}}
\put(3976,-2836){\makebox(0,0)[lb]{\smash{{\SetFigFont{5}{6.0}{\rmdefault}{\mddefault}{\updefault}{\color[rgb]{0,0,0}$B$}%
}}}}
\put(-1400,-2500){\makebox(0,0)[lb]{\smash{{\SetFigFont{5}{6.0}{\rmdefault}{\mddefault}{\updefault}{\color[rgb]{0,0,0}$a$}%
}}}}
\put(-800,-3586){\makebox(0,0)[lb]{\smash{{\SetFigFont{5}{6.0}{\rmdefault}{\mddefault}{\updefault}{\color[rgb]{0,0,0}$A$}%
}}}}
\put(3976,-4261){\makebox(0,0)[lb]{\smash{{\SetFigFont{5}{6.0}{\rmdefault}{\mddefault}{\updefault}{\color[rgb]{0,0,0}$B$}%
}}}}
\put(-5964,-4486){\makebox(0,0)[lb]{\smash{{\SetFigFont{5}{6.0}{\rmdefault}{\mddefault}{\updefault}{\color[rgb]{0,0,0}$15$}%
}}}}
\put(-5964,-2686){\makebox(0,0)[lb]{\smash{{\SetFigFont{5}{6.0}{\rmdefault}{\mddefault}{\updefault}{\color[rgb]{0,0,0}$14$}%
}}}}
\put(-4164,-2686){\makebox(0,0)[lb]{\smash{{\SetFigFont{5}{6.0}{\rmdefault}{\mddefault}{\updefault}{\color[rgb]{0,0,0}$13$}%
}}}}
\put(-2364,-2686){\makebox(0,0)[lb]{\smash{{\SetFigFont{5}{6.0}{\rmdefault}{\mddefault}{\updefault}{\color[rgb]{0,0,0}$12$}%
}}}}
\put(-2364,-4486){\makebox(0,0)[lb]{\smash{{\SetFigFont{5}{6.0}{\rmdefault}{\mddefault}{\updefault}{\color[rgb]{0,0,0}$17$}%
}}}}
\put(-6074,-3586){\makebox(0,0)[lb]{\smash{{\SetFigFont{5}{6.0}{\rmdefault}{\mddefault}{\updefault}{\color[rgb]{0,0,0}$b$}%
}}}}
\put(-5800,-3586){\makebox(0,0)[lb]{\smash{{\SetFigFont{5}{6.0}{\rmdefault}{\mddefault}{\updefault}{\color[rgb]{0,0,0}$A$}%
}}}}
\put(-5024,-2500){\makebox(0,0)[lb]{\smash{{\SetFigFont{5}{6.0}{\rmdefault}{\mddefault}{\updefault}{\color[rgb]{0,0,0}$A$}%
}}}}
\put(-3224,-2836){\makebox(0,0)[lb]{\smash{{\SetFigFont{5}{6.0}{\rmdefault}{\mddefault}{\updefault}{\color[rgb]{0,0,0}$a$}%
}}}}
\put(-3264,-2536){\makebox(0,0)[lb]{\smash{{\SetFigFont{5}{6.0}{\rmdefault}{\mddefault}{\updefault}{\color[rgb]{0,0,0}$B$}%
}}}}
\put(-1424,-4336){\makebox(0,0)[lb]{\smash{{\SetFigFont{5}{6.0}{\rmdefault}{\mddefault}{\updefault}{\color[rgb]{0,0,0}$b$}%
}}}}
\put(-1464,-4636){\makebox(0,0)[lb]{\smash{{\SetFigFont{5}{6.0}{\rmdefault}{\mddefault}{\updefault}{\color[rgb]{0,0,0}$A$}%
}}}}
\put(-564,-4486){\makebox(0,0)[lb]{\smash{{\SetFigFont{5}{6.0}{\rmdefault}{\mddefault}{\updefault}{\color[rgb]{0,0,0}$16$}%
}}}}
\put(-564,-2686){\makebox(0,0)[lb]{\smash{{\SetFigFont{5}{6.0}{\rmdefault}{\mddefault}{\updefault}{\color[rgb]{0,0,0}$11$}%
}}}}
\put(1276,-2686){\makebox(0,0)[lb]{\smash{{\SetFigFont{5}{6.0}{\rmdefault}{\mddefault}{\updefault}{\color[rgb]{0,0,0}$4$}%
}}}}
\put(1276,-4486){\makebox(0,0)[lb]{\smash{{\SetFigFont{5}{6.0}{\rmdefault}{\mddefault}{\updefault}{\color[rgb]{0,0,0}$7$}%
}}}}
\put(3076,-4486){\makebox(0,0)[lb]{\smash{{\SetFigFont{5}{6.0}{\rmdefault}{\mddefault}{\updefault}{\color[rgb]{0,0,0}$8$}%
}}}}
\put(4876,-886){\makebox(0,0)[lb]{\smash{{\SetFigFont{5}{6.0}{\rmdefault}{\mddefault}{\updefault}{\color[rgb]{0,0,0}$3$}%
}}}}
\put(3076,-2686){\makebox(0,0)[lb]{\smash{{\SetFigFont{5}{6.0}{\rmdefault}{\mddefault}{\updefault}{\color[rgb]{0,0,0}$5$}%
}}}}
\put(4876,-2686){\makebox(0,0)[lb]{\smash{{\SetFigFont{5}{6.0}{\rmdefault}{\mddefault}{\updefault}{\color[rgb]{0,0,0}$6$}%
}}}}
\put(4876,-4486){\makebox(0,0)[lb]{\smash{{\SetFigFont{5}{6.0}{\rmdefault}{\mddefault}{\updefault}{\color[rgb]{0,0,0}$9$}%
}}}}
\put(6640,-4486){\makebox(0,0)[lb]{\smash{{\SetFigFont{5}{6.0}{\rmdefault}{\mddefault}{\updefault}{\color[rgb]{0,0,0}$10$}%
}}}}
\put(2200,-4636){\makebox(0,0)[lb]{\smash{{\SetFigFont{5}{6.0}{\rmdefault}{\mddefault}{\updefault}{\color[rgb]{0,0,0}$b$}%
}}}}
\put(2140,-4336){\makebox(0,0)[lb]{\smash{{\SetFigFont{5}{6.0}{\rmdefault}{\mddefault}{\updefault}{\color[rgb]{0,0,0}$A$}%
}}}}
\put(5776,-4336){\makebox(0,0)[lb]{\smash{{\SetFigFont{5}{6.0}{\rmdefault}{\mddefault}{\updefault}{\color[rgb]{0,0,0}$a$}%
}}}}
\put(5776,-4636){\makebox(0,0)[lb]{\smash{{\SetFigFont{5}{6.0}{\rmdefault}{\mddefault}{\updefault}{\color[rgb]{0,0,0}$B$}%
}}}}
\put(4000,-2536){\makebox(0,0)[lb]{\smash{{\SetFigFont{5}{6.0}{\rmdefault}{\mddefault}{\updefault}{\color[rgb]{0,0,0}$a$}%
}}}}
\put(6151,-2986){\makebox(0,0)[lb]{\smash{{\SetFigFont{5}{6.0}{\rmdefault}{\mddefault}{\updefault}{\color[rgb]{0,0,0}$b$}%
}}}}
\put(4726,-1861){\makebox(0,0)[lb]{\smash{{\SetFigFont{5}{6.0}{\rmdefault}{\mddefault}{\updefault}{\color[rgb]{0,0,0}$b$}%
}}}}
\put(3751,-1711){\makebox(0,0)[lb]{\smash{{\SetFigFont{5}{6.0}{\rmdefault}{\mddefault}{\updefault}{\color[rgb]{0,0,0}$A$}%
}}}}
\put(1426,-1786){\makebox(0,0)[lb]{\smash{{\SetFigFont{5}{6.0}{\rmdefault}{\mddefault}{\updefault}{\color[rgb]{0,0,0}$A$}%
}}}}
\put(2140,-2500){\makebox(0,0)[lb]{\smash{{\SetFigFont{5}{6.0}{\rmdefault}{\mddefault}{\updefault}{\color[rgb]{0,0,0}$B$}%
}}}}
\put(1051,-3586){\makebox(0,0)[lb]{\smash{{\SetFigFont{5}{6.0}{\rmdefault}{\mddefault}{\updefault}{\color[rgb]{0,0,0}$b$}%
}}}}
\end{picture}%

%% file: jsc-3474-fig5.pstex_t
\begin{picture}(0,0)%
\includegraphics{jsc-3474-fig5.pstex}%
\end{picture}%
\setlength{\unitlength}{1579sp}%
\begingroup\makeatletter\ifx\SetFigFont\undefined%
\gdef\SetFigFont#1#2#3#4#5{%
  \reset@font\fontsize{#1}{#2pt}%
  \fontfamily{#3}\fontseries{#4}\fontshape{#5}%
  \selectfont}%
\fi\endgroup%
\begin{picture}(12916,7863)(-2707,-8383)
\put(6076,-736){\makebox(0,0)[lb]{\smash{{\SetFigFont{8}{9.6}{\familydefault}{\mddefault}{\updefault}{\color[rgb]{0,0,0}$Y(\bN)$}%
}}}}
\put(-2249,-6736){\makebox(0,0)[lb]{\smash{{\SetFigFont{8}{9.6}{\familydefault}{\mddefault}{\updefault}{\color[rgb]{0,0,0}$c$}%
}}}}
\put(-1274,-8011){\makebox(0,0)[lb]{\smash{{\SetFigFont{8}{9.6}{\familydefault}{\mddefault}{\updefault}{\color[rgb]{0,0,0}$\bM$}%
}}}}
\put(1726,-6211){\makebox(0,0)[lb]{\smash{{\SetFigFont{8}{9.6}{\familydefault}{\mddefault}{\updefault}{\color[rgb]{0,0,0}$F$}%
}}}}
\put(3376,-6136){\makebox(0,0)[lb]{\smash{{\SetFigFont{8}{9.6}{\familydefault}{\mddefault}{\updefault}{\color[rgb]{0,0,0}$F(c)$}%
}}}}
\put(5176,-5611){\makebox(0,0)[lb]{\smash{{\SetFigFont{8}{9.6}{\familydefault}{\mddefault}{\updefault}{\color[rgb]{0,0,0}$F(m)$}%
}}}}
\put(6451,-5311){\makebox(0,0)[lb]{\smash{{\SetFigFont{8}{9.6}{\familydefault}{\mddefault}{\updefault}{\color[rgb]{0,0,0}$F(c')$}%
}}}}
\put(7951,-5686){\makebox(0,0)[lb]{\smash{{\SetFigFont{8}{9.6}{\familydefault}{\mddefault}{\updefault}{\color[rgb]{0,0,0}$nn'$}%
}}}}
\put(9151,-6136){\makebox(0,0)[lb]{\smash{{\SetFigFont{8}{9.6}{\familydefault}{\mddefault}{\updefault}{\color[rgb]{0,0,0}$b'$}%
}}}}
\put(7576,-7636){\makebox(0,0)[lb]{\smash{{\SetFigFont{8}{9.6}{\familydefault}{\mddefault}{\updefault}{\color[rgb]{0,0,0}$b$}%
}}}}
\put(5101,-6961){\makebox(0,0)[lb]{\smash{{\SetFigFont{8}{9.6}{\familydefault}{\mddefault}{\updefault}{\color[rgb]{0,0,0}$F(m)n$}%
}}}}
\put(8476,-6886){\makebox(0,0)[lb]{\smash{{\SetFigFont{8}{9.6}{\familydefault}{\mddefault}{\updefault}{\color[rgb]{0,0,0}$n'$}%
}}}}
\put(6226,-8311){\makebox(0,0)[lb]{\smash{{\SetFigFont{8}{9.6}{\familydefault}{\mddefault}{\updefault}{\color[rgb]{0,0,0}$\bN$}%
}}}}
\put(-1424,-886){\makebox(0,0)[lb]{\smash{{\SetFigFont{8}{9.6}{\familydefault}{\mddefault}{\updefault}{\color[rgb]{0,0,0}$X(\bM)$}%
}}}}
\put(-1499,-4561){\makebox(0,0)[lb]{\smash{{\SetFigFont{8}{9.6}{\familydefault}{\mddefault}{\updefault}{\color[rgb]{0,0,0}$X$}%
}}}}
\put(-2174,-2836){\makebox(0,0)[lb]{\smash{{\SetFigFont{8}{9.6}{\familydefault}{\mddefault}{\updefault}{\color[rgb]{0,0,0}$X(c)$}%
}}}}
\put(-224,-2236){\makebox(0,0)[lb]{\smash{{\SetFigFont{8}{9.6}{\familydefault}{\mddefault}{\updefault}{\color[rgb]{0,0,0}$X(c')$}%
}}}}
\put(-974,-2836){\makebox(0,0)[lb]{\smash{{\SetFigFont{8}{9.6}{\familydefault}{\mddefault}{\updefault}{\color[rgb]{0,0,0}$X(m)$}%
}}}}
\put(-1349,-5986){\makebox(0,0)[lb]{\smash{{\SetFigFont{8}{9.6}{\familydefault}{\mddefault}{\updefault}{\color[rgb]{0,0,0}$m$}%
}}}}
\put(-149,-6286){\makebox(0,0)[lb]{\smash{{\SetFigFont{8}{9.6}{\familydefault}{\mddefault}{\updefault}{\color[rgb]{0,0,0}$c'$}%
}}}}
\put(6976,-6436){\makebox(0,0)[lb]{\smash{{\SetFigFont{8}{9.6}{\familydefault}{\mddefault}{\updefault}{\color[rgb]{0,0,0}$n$}%
}}}}
\put(5626,-1636){\makebox(0,0)[lb]{\smash{{\SetFigFont{8}{9.6}{\familydefault}{\mddefault}{\updefault}{\color[rgb]{0,0,0}$Y(F(c'))$}%
}}}}
\put(6001,-2611){\makebox(0,0)[lb]{\smash{{\SetFigFont{8}{9.6}{\familydefault}{\mddefault}{\updefault}{\color[rgb]{0,0,0}$Y(n)$}%
}}}}
\put(6901,-3436){\makebox(0,0)[lb]{\smash{{\SetFigFont{8}{9.6}{\familydefault}{\mddefault}{\updefault}{\color[rgb]{0,0,0}$Y(b)$}%
}}}}
\put(6001,-4486){\makebox(0,0)[lb]{\smash{{\SetFigFont{8}{9.6}{\familydefault}{\mddefault}{\updefault}{\color[rgb]{0,0,0}$Y$}%
}}}}
\end{picture}%